\documentclass[11pt]{article} 
\usepackage{amsmath}
\usepackage{amssymb,amsmath, amsfonts, amsthm}
\usepackage{graphicx}
\usepackage[shortlabels]{enumitem}
\usepackage{algorithm}
\usepackage{algorithmic}
\usepackage{tikz}
\usepackage{url}
\usepackage[makeroom]{cancel}

\usepackage{cancel}
\usepackage[margin=1.5in]{geometry}
\usepackage{mathtools}
\usepackage[font=footnotesize]{caption}

\usepackage{url}
\usepackage{commath}
\newtheorem*{theorem*}{Theorem}

\tikzstyle{vertex}=[circle, draw, inner sep=0pt, minimum size=6pt]
\usetikzlibrary{shapes.geometric, arrows}
\usetikzlibrary{graphs,graphs.standard}
\tikzset{main node/.style={circle,fill=blue!20,draw,minimum size=.5cm,inner sep=0pt},
            }
\newcommand{\R}{\ensuremath{\mathbb{R}}}

\usepackage{verbatim}  
\usepackage{graphicx}  
\usepackage{float}         
\usepackage{titling}
\usepackage{blindtext}
\usepackage{multirow}
\usepackage{esvect}
\usepackage{multirow}
\usepackage{bm}
\usepackage{afterpage}

\usepackage[symbol]{footmisc}

\renewcommand{\thefootnote}{\fnsymbol{footnote}}

\numberwithin{equation}{section}

 \usepackage{setspace}
\title{An Assessment of the Supremizer and Aggregation Methods of Stabilization for Reduced-Order Models}
\author{Kayla D. Davie$^\ast$  and Howard C. Elman$^\dagger$}

\date{}

\begin{document}

\maketitle

\footnotetext[1]{Applied Mathematics \& Statistics, and Scientific Computation Program, University of Maryland, College Park, MD 20742}
\footnotetext[2]{Department of Computer Science and Institute for Advanced Computer Studies, University of Maryland, College Park, MD\ 20742}
\renewcommand*{\thefootnote}{\arabic{footnote}}

\vspace{-.6in}

\abstract{
\begin{footnotesize}
We explore the features of two methods of stabilization, aggregation and supremizer methods, for reduced-order modeling of parametrized optimal control problems. In both methods, the reduced basis spaces are augmented to guarantee stability. For the aggregation method, the reduced basis approximation spaces for the state and adjoint variables are augmented in such a way that the spaces are identical. For the supremizer method, the reduced basis approximation space for the state-control product space is augmented with the solutions of a supremizer equation. We implement both of these methods for solving several parametrized control problems and assess their performance. Results indicate that the number of reduced basis vectors needed to approximate the solution space to some tolerance with the supremizer method is much larger, possibly double, that for aggregation. There are also some cases where the supremizer method fails to produce a converged solution.
We present results to compare the accuracy, efficiency, and computational costs associated with both methods of stabilization which suggest that stabilization by aggregation is a superior stabilization method for control problems.
\end{footnotesize}
}

\vspace{.1in} 


\section{Introduction}

Consider the parametrized control problem subject to a constraint given by a partial differential equation (PDE)
\begin{align}\label{eq:genformppdecontrol}
    \min_{u, f} ~~ \mathcal{J}(u, f; \mu) ~~~~~~ \text{subject to} ~~~~~~ \mathcal{G}(u, f; \mu)=0.
\end{align}
Here, $\mathcal{J}$ is the cost functional to be minimized, $u$ is the state variable, $f$ is the control variable, $\mu$ is a vector of parameters, and $\mathcal{G}$ is the PDE constraint including boundary conditions. 
The algebraic system that results from the discretization of (\ref{eq:genformppdecontrol}) 
is of large dimension and saddle-point form. Iterative Krylov subspace methods have been proposed for efficiently approximating the solution to deterministic versions of the problem \cite{esw, haberascher, hein_lq, laz, rdw}.
In the parametrized setting, solutions to these problems are often required for a large number of parameters. This may be because a simulation is done to explore solutions for a large number of parameters, or to obtain statistical properties of solutions. Thus, this computationally expensive problem must be solved repeatedly, a major computational task.

Reduced order modeling (ROM) is an efficient way to address this issue \cite{hesthaven, quarteronimanzoninegri}. This technique uses a reduced basis (RB) method to replace a problem of high dimension with one of reduced order without jeopardizing accuracy of the approximation. These reduced-order models require less storage and are significantly less computationally taxing than the ``full-sized" version of the models. RB methods typically consist of two stages, referred to as an \textit{offline} stage and an \textit{online} stage. In the offline stage, the reduced basis is constructed. This stage requires obtaining a set of ``snapshots" or ``truth" approximations of the parametrized optimal control problem for various choices of parameters, consisting of high-dimensional discretizations of the solution that can be acquired using standard discretization methods such as the finite element method. This stage may be computationally expensive, but once the reduced basis is available, the online stage consists of simulation for multiple values of the parameter, which can be done at lower costs, so that the overhead of the online computations is amortized over these computations. 

For constrained problems such as (\ref{eq:genformppdecontrol}), 
standard approaches for constructing reduced-order models require special treatment of the basis to ensure inf-sup stability of the reduced model. Two common ways of augmenting the reduced basis to satisfy the condition are through use of a supremizer \cite{quarteronirozza, rozzaveroy}
and through aggregation \cite{nrmq}. 
Both approaches have been studied in depth and have been shown to produce reduced bases that are reliable. 
The goal of this paper is to study both of these methods of stabilization and compare the accuracy, efficiency, and computational costs associated with them. We look at two benchmark problems, parametrized versions of the diffusion control problem and the convection-diffusion control problem.
Our main observation is that stabilization by aggregation results in a basis that requires significantly fewer snapshots with less work required to use the basis in the online stage.

An outline of this paper is as follows. In Section \ref{sec:paroptcont}, we will define the operators, matrix system, and spaces associated with the deterministic diffusion control problem, and we review the inf-sup condition for the deterministic problem. In Section \ref{sec:rom}, we discuss reduced-order modeling and describe the algorithm used to construct a reduced basis. In Section \ref{sec:stabilization}, we present the two different methods of stabilization and review how inf-sup stability is established for both RB systems. In Section \ref{sec:numres}, we present our numerical results for the two benchmark problems.

\section{The Parametrized Optimal Control Problem}\label{sec:paroptcont}

Consider the parametrized diffusion control problem,
\begin{align}
\min_{u,f} ~ \frac{1}{2} \norm{u(x, \mu) - \hat{u}(x, \mu)}_{L_{2}(\Omega)}^{2} + \beta \norm{f(x,\mu)}_{L_{2}(\Omega)}^{2} \label{eq:pardiffctrl1}\\
\text{subject to} \quad -\bigtriangledown \cdot ( \sigma(x, \mu) \bigtriangledown u(x, \mu) )~ = ~ f(x, \mu) ~ \text{ in }  ~ \Omega \times \Gamma, \label{eq:pardiffctrl2}\\
\quad \quad \text{such that} \quad \quad \quad  u(x, \mu)  = g(x) \quad  ~ \text{on} ~ \partial \Omega_{D} \times \Gamma,  \quad \quad    \quad \quad \quad \quad \label{eq:pardiffctrl3}\\
\quad \quad \quad \quad \quad  \sigma(x, \mu) \frac{\partial  u(x, \mu)}{\partial n}  = 0 \quad  ~ \text{on} ~ \partial \Omega_{N} \times \Gamma.    \quad \quad \quad \quad \quad \quad \label{eq:pardiffctrl4}
\end{align}
Here, $\Omega \subset \R^d$ is the spatial domain, $\mu$ is a vector of parameters, and $\sigma(x, \mu)$ is a parameter-dependent diffusion coefficient.
For each $\mu$, $\hat{u}(x, \mu)$ is a desired state and we seek a state, $u(x, \mu)$, as close to $\hat{u}(x, \mu)$ as possible (in the $L_2$-norm sense), and the control, $f(x, \mu)$, can vary to achieve this. The  regularization term (second term in cost functional) is added to make the problem well-posed \cite{haberascher, rdw}.  
This term involves a regularization parameter, $\beta$, based on Tikhonov regularization theory \cite{vondaniels},
and is typically set to be $\beta \approx 10^{-2}$. The solution variables are state and control variables, $u(x, \mu)$ and $f(x, \mu)$ respectively, where $\hat{u}(x; \mu)$ and $g(x)$ are given.
We will consider a two-dimensional spatial domain $\Omega$ divided into $N_D$ horizontal strips or subdomains with the $k$th subdomain denoted $\Omega_k$.
The diffusion coefficient is taken to be piecewise constant on each subdomain $\sigma(x,\mu) |_{\Omega_k} = \mu_k, k=1:N_D,$
giving a parameter space of dimension $N_D$ depending on the $N_D$-dimensional parameter vector $\mu = [\mu_1,..., \mu_{N_D}]^T$.

\subsection{Definition of Operators and Spaces}\label{sec:defineops}

In this section, we introduce the operators and spaces for the parameterized diffusion control problem with fixed $\mu$ and present the structure for solving a problem of this type at the continuous level. The goal is to find state and control solutions, $u(x, \mu)$ and $f(x, \mu)$ respectively, for some desired state $\hat{u}(x, \mu)$ and Dirichlet boundary data $g(x)$. 
The state space, $U$, is a Hilbert space for state function $u(x, \mu)$ equipped with inner product and induced norm 
\begin{align}\label{eq:sharednorm}
    (u_1, u_2)_U = \int_{\Omega} \bigtriangledown u_1(x, \mu) \cdot \bigtriangledown u_2(x, \mu) dx,~~ ||u||=(
    u, u )_U^{1/2}, \quad
\end{align}
the $H^1$ seminorm. The solution space for $u(x, \mu)$, $H^1_E$, has the Dirichlet condition built into its definition, and functions in the test space $H^1_{E_0}$ are zero on the Dirichlet part of boundary $\Omega_D$. 
The control space, $F$, for control function $f(x, \mu)$, is another Hilbert space equipped with inner product and induced norm 
\begin{align*}
    (f_1, f_2)_F = \int_{\Omega}  f_1(x, \mu) f_2(x, \mu) dx, ~~ ||f||=(
    f, f )_F^{1/2}, \quad
\end{align*}
i.e. $F=L_2(\Omega)$. A standard approach to solving constrained optimization problems is to apply first-order conditions for stationarity to the Lagrangian functional
\begin{align}\label{eq:lagrangian}
    \mathrm{L} \! := \! \frac{1}{2} \int_{\Omega} \! (u(x, \mu) \! - \! \hat{u}(x, \mu))^{2} dx + \beta \! \int_{\Omega} \! f(x, \mu)^{2} dx + \lambda \Big(\int_{\Omega} \! \sigma(x, \mu) \bigtriangledown u(x, \mu) \! \cdot \! \bigtriangledown v(x, \mu) dx -  \! \int_{\Omega} \!  f(x, \mu) v(x, \mu) dx \Big).
\end{align}
This gives rise to an adjoint/Lagrange multiplier function, $\lambda(x, \mu)$, and the adjoint/Lagrange multiplier space, $Q$. 
The Lagrange multiplier satisfies \cite{esw} 
$$- \bigtriangledown^2 \lambda(x, \mu) = -u(x, \mu) + \hat{u}(x, \mu), $$ 
 with homogeneous Dirichlet boundary conditions,
 so that  $\lambda(x, \mu) \in H_{E_0}^1$  as well. 
 Thus, the Hilbert space $Q$ is equipped with inner product and induced norm (\ref{eq:sharednorm}) 
 and $Q$ and $U$ are equivalent.

It will be convenient to refer to certain product spaces for the purpose of defining operators and establishing stability. These include the product between control and state spaces, $X=F \times U$, with elements $\bar{x}=(f(x, \mu), u(x, \mu)) \in X$, 
and the product spaces $U \times Q$, $F \times Q$, and $X \times Q$. For given $\mu$, let $a( \cdot, \cdot; \mu):  U \times Q \rightarrow \R$ be the linear elliptic operator 
$$a(z, q; \mu) = \int_{\Omega} \sigma(x, \mu) \bigtriangledown z(x, \mu) \cdot \bigtriangledown q(x, \mu) dx.$$ 
We assume that the operator $a$ is coercive, i.e., there exists $ \alpha_0>0$ such that $\displaystyle \alpha(\mu) = inf_{z \in U} \frac{a(z,z;\mu)}{||z||_U^2} \geq \alpha_0$ for all $\mu$.
The action of the control $f(x, \mu)$ is represented by the operator $c(\cdot, \cdot; \mu): F \times Q \rightarrow \R$,  
$$c(f,q; \mu) = \int_{\Omega} f(x,\mu)  q(x, \mu) dx.$$

The weak formulation of (\ref{eq:pardiffctrl1})-(\ref{eq:pardiffctrl4}) is to find the minimizers $u( \cdot, \mu) \in H^1_E(\Omega)$,  $f(\cdot, \mu) \in F$ as in (\ref{eq:pardiffctrl1}) satisfying (\ref{eq:pardiffctrl3})-(\ref{eq:pardiffctrl4}) subject to the weak form of the constraint 
 \begin{align}\label{eq:weakconstraint}
  \int_{\Omega} \sigma(x, \mu) \bigtriangledown u(x, \mu) \cdot \bigtriangledown v(x, \mu) dx = \int_{\Omega}  f(x, \mu) v(x, \mu) dx ~~ \forall v \in H^1_{E_0}(\Omega).    
 \end{align}
First-order stationarity of the Lagrangian (\ref{eq:lagrangian}) results in a saddle-point system.

It is well known that to guarantee the existence, uniqueness and stability of solutions, saddle-point systems of this form must satisfy an inf-sup condition. 
 Specifically, for the bilinear form $\mathrm{B}(\cdot, \cdot; \mu): X \times Q \rightarrow \R$, given by $$\mathrm{B}(\bar{w}, q; \mu) = a(z, q; \mu)- c(v, q; \mu)$$ where $\bar{w}={\footnotesize\begin{bmatrix} v \\ z \end{bmatrix}} \in X=F \times U$, it must hold that there exists there exists $\beta_0>0$ such that $$\text{inf}_{q \in Q} ~ \text{sup}_{\bar{w} \in X} ~ \frac{\mathrm{B}(\bar{w}, q; \mu)}{||\bar{w}||_X ||q||_Q} \geq \beta_0.$$ It is shown in \cite{nrmq} that this condition holds for the control problem (\ref{eq:pardiffctrl1})-(\ref{eq:pardiffctrl4}).

\subsection{Discrete Forms and Matrix Systems}\label{sec:dermat}

For the discrete version of the parametrized control problem, let $S_0^h$ be a finite-dimensional subspace of $H_{E_0}^h$, and let $S_E^h$ augment $S_0^h$ with a finite set of basis functions used to impose Dirichlet boundary conditions. 
 The Galerkin finite element method is chosen for discretization of the optimal control problem. The weak formulation of the constraint is then to find $u_h(x, \mu) \in S_E^h $ such that 
\begin{align}\label{eq:discreteconstraint}
\int_{\Omega} \sigma(x; \mu) \bigtriangledown u_h(x, \mu) \cdot \bigtriangledown v_h(x, \mu) = \int_{\Omega} f_h(x, \mu) v_h(x, \mu) ~~ \forall v_h \in S_0^h \subset H^1_{E_0}(\Omega).
\end{align}
Here, $S_E^h$ is the solution space and $S_{0}^{h}$ is a vector space containing test functions. Let the basis for the test functions be denoted $\{ \phi_{1},...,\phi_{n} \}$, and assume this basis is extended by $\phi_{n+1},...,\phi_{n+\partial n}$, which ensures that the Dirichlet boundary conditions hold at certain points on $\partial \Omega_D$; see \cite[pp.~ 195ff.]{brennerscott} for additional discussion of this. The discrete analog of the weak version of the minimization problem (\ref{eq:pardiffctrl1})-(\ref{eq:pardiffctrl4}) is as follows: 
\begin{align}
    \displaystyle \min_{u_h,f_h} ~& \frac{1}{2} \begin{Vmatrix} u_h(x, \mu) - \hat{u}_h(x, \mu) \end{Vmatrix}_{2}^{2} + \beta \begin{Vmatrix} f_h(x, \mu) \end{Vmatrix}_{2}^{2} \\
    ~\text{subject to}~ &\int_{\Omega} \sigma(x; \mu) \bigtriangledown u_h(x, \mu) \cdot \bigtriangledown v_h(x, \mu) dx = \int_{\Omega}  f_h(x, \mu) v_h(x, \mu) dx
\end{align}
subject to boundary conditions for $u_h$ as in (\ref{eq:pardiffctrl3})-(\ref{eq:pardiffctrl4}).

For discretization, we use equal-order $Q_1$ (piecewise bilinear) finite elements for all three (state, control and Lagrange multiplier) spaces, which is shown to be div-stable in \cite{nrmq}. In the state space, the finite element approximation of $u(x, \mu)$, $u_h(x, \mu)$, is represented in terms of basis functions $u_{h}(x, \mu)= \sum_{j=1}^{n} \textbf{u}_{\mu, j} \phi_{j} + \sum_{j=n+1}^{n+ \partial n} \textbf{u}_{\mu, j} \phi_{j} $ such that $u_{h}(x, \mu)$ is determined uniquely by coefficient vector $\textbf{u}(\mu)=\textbf{u}_{\mu}=  (\textbf{u}_{ \mu, 1},...,\textbf{u}_{\mu, n + \partial n})^{T}$.
The coefficient vector $\textbf{f}(\mu)=\textbf{f}_{\mu} = (\textbf{f}_{\mu, 1},...\textbf{f}_{\mu, n})^{T}$ associated with the approximation  $f_{h}(x, \mu)$ of $f(x, \mu)$ is defined similarly.
The cost functional can be rewritten in matrix notation as 
$$\min_{\textbf{u}_{\mu},\textbf{f}_{\mu}} ~ \frac{1}{2} \textbf{u}_{\mu}^T \mathcal{M} \textbf{u}_{\mu} - \textbf{u}_{\mu}^T \textbf{b}_{\mu} + \beta \textbf{f}_{\mu}^{T} \mathcal{M} \textbf{f}_{\mu}$$ 
where $\textbf{f}(\mu)=\textbf{f}_{\mu} = (\textbf{f}_{\mu, 1},...\textbf{f}_{\mu, n})^{T}$, $\textbf{b}(\mu)=\textbf{b}_{\mu} = \{ \int \hat{u}(x, \mu)  \phi_{i} \}_{i=1,...,n}$, and the mass matrix, $\mathcal{M}$, is defined as $\mathcal{M} = \{ \int \phi_{i} \phi_{j} \} _{i,j=1,...,n}.$ The constraint can be expressed as 
\begin{equation}\label{eq:discconstraintstokes}
    \mathcal{K}(\mu) \textbf{u}_{\mu} = \mathcal{M} \textbf{f}_{\mu} + \textbf{d}_{\mu}
\end{equation}
such that $\mathcal{K}(\mu) = [k^{\mu}_{ij}]$ and $[k^{\mu}_{ij}]=  \int_{\Omega} \sigma(x, \mu)  \bigtriangledown \phi_j \cdot \bigtriangledown \phi_i = \sum_q \mu_q \int_{\Omega_q} \bigtriangledown \phi_j \cdot \bigtriangledown \phi_i$. Note that $\textbf{d}(\mu)=\textbf{d}_{\mu}$ contains terms from the boundary values of the discretized $u(x, \mu)$, such that $\textbf{d}_{\mu} = \{ \textbf{d}_{\mu, i} \}_{i=1,...,n}$ with $\textbf{d}_{\mu, i}= \int_{\Omega} g \phi_{i} ~d\Omega - \sum_{j=n+1}^{n+n_{\partial}} \textbf{u}_{\mu, j} \int_{\Omega} \sigma(x, \mu) \bigtriangledown \phi_{i} \cdot \bigtriangledown \phi_{j} ~ d \Omega$. 
The discrete Lagrangian is then
$$\mathrm{L} := \frac{1}{2} \textbf{u}_{\mu}^{T} \mathcal{M} \textbf{u}_{\mu} - \textbf{u}_{\mu}^{T} \textbf{b}_{\mu} + \beta \textbf{f}_{\mu}^{T} \mathcal{M} \textbf{f}_{\mu} + \boldsymbol{\lambda}_{\mu}^{T} (\mathcal{K}(\mu) \textbf{u}_{\mu} - \mathcal{M} \textbf{f}_{\mu} - \textbf{d}_{\mu})$$
where  $\boldsymbol{\lambda}_{\mu}$ is a vector of Lagrange multipliers associated with finite element approximation of $\lambda(x, \mu)$, $\lambda_h(x, \mu)$.

Applying first-order conditions for stationarity 
gives a set of three coupled equations, a block system 
\begin{equation}\label{eq:poisssys}
\begin{bmatrix} 2 \beta \mathcal{M} & 0 & -\mathcal{M} \\ 0 & \mathcal{M} & \mathcal{K}(\mu)^{T} \\ -\mathcal{M} & \mathcal{K}(\mu) & 0 \end{bmatrix} \begin{bmatrix} \textbf{f}_{\mu} \\ \textbf{u}_{\mu} \\ \boldsymbol{\lambda}_{\mu} \end{bmatrix} =  \begin{bmatrix} 0 \\ \textbf{b}_{\mu} \\ \textbf{d}_{\mu} \end{bmatrix}.
\end{equation}
The block $3 \times 3$ coefficient matrix $\mathcal{G}(\mu)$ can also be represented in compact form,
\begin{equation}\label{eq:saddle}
   \mathcal{G}(\mu)= \begin{bmatrix}
    \mathcal{A} & \mathcal{B}(\mu)^T \\ \mathcal{B}(\mu) & 0
    \end{bmatrix},
\end{equation}
where $ \mathcal{B}(\mu)=[-\mathcal{M}, \mathcal{K}(\mu)]$, and the matrix system can be denoted $\mathcal{G}(\mu) \mathbf{v}(\mu)=\mathbf{b}(\mu)$.

The discrete spaces are defined as follows. The state space, $U_h$, has inner product $\langle u_h(x, \mu), v_h(x, \mu) \rangle_{U_h} =  \langle \mathcal{K} \textbf{u}_{\mu}, \textbf{v}_{\mu} \rangle = \textbf{v}_{\mu}^T \mathcal{K} \textbf{u}_{\mu}$ and induced norm $|| u_h||_{U_h}=\langle \mathcal{K} \textbf{u}_{\mu}, \textbf{u}_{\mu} \rangle^{1/2}$  where $\mathcal{K}$ is defined like $\mathcal{K}(\mu)$ above with $\sigma \equiv 1$. 
The control space, $F_h$, is equipped with inner product $\langle f_h(x, \mu), v_h(x, \mu) \rangle_{F_h} = \langle \mathcal{M} \textbf{f}_{\mu}, \textbf{v}_{\mu} \rangle = \textbf{v}_{\mu}^T \mathcal{M} \textbf{f}_{\mu}$ and induced norm $|| f_h(x, \mu)||_{F_h}=\langle  \mathcal{M} \textbf{f}_{\mu}, \textbf{f}_{\mu} \rangle^{1/2}.$
The Lagrange multiplier space $Q_h$  is equipped with inner product $\langle q_h(x, \mu), p_h(x, \mu) \rangle_{Q_h}= \langle \mathcal{K} \textbf{q}_{\mu}, \textbf{p}_{\mu} \rangle = \textbf{p}_{\mu}^T \mathcal{K} \textbf{q}_{\mu}$ and induced norm $ || \lambda_h(x, \mu)||_{Q_h}= \langle \mathcal{K} \boldsymbol{\lambda}_{\mu}, \boldsymbol{\lambda}_{\mu} \rangle^{1/2}.$ The state-control space, $X_h=F_h \times U_h$, has elements $\bar{x}_h$ with associated coefficient vector $\bar{\textbf{x}}(\mu)={\footnotesize\begin{bmatrix}
\textbf{f}_{\mu}\\ \textbf{u}_{\mu}
\end{bmatrix}}$.

Because the matrix system (\ref{eq:poisssys}) is of saddle point form (\ref{eq:saddle}), satisfaction of the inf-sup condition is required to guarantee stability. The discrete inf-sup condition requires that there exists $\beta_0>0$ such that $$\min_{\textbf{q} \in Q_h} ~ \max_{\bar{\textbf{w}}={\tiny\begin{bmatrix}
    \textbf{v} \\ \textbf{z}
\end{bmatrix}} \in X_h} ~ \frac{\langle \mathcal{K}(\mu) \textbf{z}, \textbf{q} \rangle - \langle \mathcal{M} \textbf{v}, \textbf{q} \rangle}{||\bar{\textbf{w}}||_X ||\textbf{q}||_Q}  \geq \beta_0,$$
as proven in \cite{nrmq}. Recall also \cite{esw} that this stability bound is equivalent to the bound on the Rayleigh quotient
\begin{align}\label{eq:schurversion}
    \frac{\langle \mathcal{B}(\mu) \mathcal{A}^{-1} \mathcal{B}(\mu)^T \mathbf{q}, \mathbf{q} \rangle^{1/2}}{\langle \mathcal{K} \mathbf{q}, \mathbf{q} \rangle^{1/2}} \geq \beta_0.
\end{align}

\section{Reduced-Order Modeling}\label{sec:rom}

ROM uses a relatively small set of full-order discrete solutions $(\textbf{f}_{\mu}, \textbf{u}_{\mu})$, called ``snapshots," 
for various values of the parameter(s) $\mu$, to approximate  the solution to the problem for other values of the parameter by projecting on the space spanned by a subset of these snapshots, 
where it is assumed that the solution manifold can be approximated well by a subset of snapshots \cite{itoravindran, nrmq, rhm}. Computations are separated into an ``offline" step in which the reduced basis is constructed, and an ``online" step in which simulations are done. We outline this process here, with additional details given in Section \ref{sec:stabilization}.

In the offline step, full solutions $(\textbf{f}_{\mu}, \textbf{u}_{\mu})$ are computed for multiple values of the parameter vector $\mu$, and these solutions are used to generate the reduced basis approximation spaces using a greedy algorithm \cite{greedyref1, greedyref2}. Let $T$ be a  set of $N_{max}$ \textit{training parameters}, where we assume that $N_{max}$ is large enough so that the span of the solution set $\Biggl\{ \footnotesize{ \begin{bmatrix} \textbf{f}_{\mu} \\ \textbf{u}_{\mu} \end{bmatrix}} \Big| \mu \in T \Biggr\}$ represents a good approximation of the entire solution space. Starting with a single parameter from $T$ and corresponding snapshot, an initial reduced basis approximation space consists of the single snapshot. At each step of the greedy algorithm for building the reduced basis, an approximation to the full solution, $\textbf{v}_r(\mu) \approx \textbf{v}(\mu)$, from the reduced space is computed for each $\mu$ in the training set, and for the parameter $\mu$ for which an error indicator $\eta(\mathbf{v}_r(\mu))$ is maximal, the snapshot, full solution $\textbf{v}(\mu)$, is added to the reduced basis. This continues until the values of the error indicator for all approximations from the reduced space are less than some prescribed tolerance for all parameters in $T$.
The resulting reduced basis approximation space is spanned by $N$ snapshots corresponding to some parameters $\mu^{(n)}, 1 \leq n \leq N$ where $N \leq N_{max}$. A generic statement of the greedy algorithm, for finding a subspace $V_i$ of $V_h \equiv F_h \times U_h \times Q_h$, is given in Algorithm \ref{alg:greedy}. At each stage, approximate solutions $\mathbf{v}_r(\mu)$ in the (current) reduced space $V_{i-1}$ are computed for all parameters $\mu$ in $T$, together with the error indicator $\eta(\mathbf{v}_r(\mu))$. For the parameter with the largest value of $\eta(\mathbf{v}_r(\mu))$, $\mu^{(i)}$, the full-system solution $\mathbf{v}(\mu^{(i)})$ is computed and the reduced basis is augmented with this solution.

\begin{algorithm}[H]
 \caption{Greedy sampling algorithm for constructing the reduced basis space $V_N$ }
 \label{alg:greedy}
 \begin{algorithmic}\label{alg:randselec}
 \STATE Randomly sample $N_{max}$ parameters $\mu$. Let $T=\{ \mu^{(k)} \}_1^{N_{max}}$.
 \STATE Choose $\mu^{(1)}$. Solve the full system for $\mathbf{v}(\mu^{1})$.
 \STATE Set $V_1=\{\mathbf{v}(\mu^{1})\}$.
 \STATE Set $N=1$.
 \STATE Set $\eta^{*}=\infty$.
 \WHILE{$\eta^{*} >$ tolerance}
 \FOR{$i=1:N_{max}$}
 \STATE Solve reduced system for $\tilde{\mathbf{v}}_r(\mu^{(i)})$, compute $\mathbf{v}_r(\mu^{(i)})=\mathcal{Q} \tilde{\mathbf{v}}_r(\mu^{(i)})$.
 \STATE Compute $\eta_i \equiv \eta(\mathbf{v}_r(\mu^{(i)}))$.
 \ENDFOR
 \STATE Find $\eta^{*}=\max_{i} \eta_i,  \mu^{*}= \arg \max_{i} \eta(\mathbf{v}_r(\mu^{(i)}))$.
 \IF{$\eta^{*} >$ tolerance}
 \STATE Set $N=N+1$.
 \STATE Solve the full system for $\mathbf{v}(\mu^{*})$. 
 \STATE Update $V_N= \text{span} (V_{N-1} \bigcup \{ \mathbf{v}(\mu^{*})\}).$ 
 \ENDIF
 \ENDWHILE 
 \end{algorithmic}
 \end{algorithm}

We comment on several things that need further discussion. First, 
we will consider two ways to specify the reduced system required for Algorithm \ref{alg:greedy}:
\begin{align}
    &\text{Galerkin projection: } [\mathcal{Q}^T \mathcal{G} (\mu) \mathcal{Q}] \tilde{\mathbf{v}}_r(\mu) = \mathcal{Q}^T \mathbf{b}(\mu) \label{eq:galsys}\\
    & \text{Petrov-Galerkin projection: } [( \mathcal{G} (\mu) \mathcal{Q})^T ( \mathcal{G} (\mu)  \mathcal{Q})] \tilde{\mathbf{v}}_r(\mu) =  ( \mathcal{G} (\mu) \mathcal{Q})^T \mathbf{b}(\mu) \label{eq:petgalsys}
\end{align}
We will discuss the details and impact of these choices in Section \ref{sec:numres}.
Second, as given, Algorithm \ref{alg:greedy} does not address the issue of inf-sup stability, which is not automatically satisfied using only snapshots. We will elaborate on this in Section \ref{sec:stabilization}.

The basis for the reduced space will ultimately be represented by the columns of a matrix $\mathcal{Q}$ such that $\mathcal{Q}\tilde{\mathbf{v}}_r=\mathbf{v}_r \approx \mathbf{v}.$ $\mathcal{Q}$ has block-diagonal structure consisting of components for the state, control and adjoint variables:
\begin{align}\label{eq:tworedmat}
\mathcal{Q} = \begin{bmatrix}
    \mathcal{Q}_f & & \\ & \mathcal{Q}_u & \\ & & \mathcal{Q}_{\lambda} 
    \end{bmatrix} ~~~ \text{or} ~~~  \mathcal{Q} = \begin{bmatrix}
    \mathcal{Q}_{\bar{x}} &  \\ & \mathcal{Q}_{\lambda}
    \end{bmatrix}    
\end{align}
The Galerkin formulation then results in one of the following scenarios:
    {\begin{align}\label{eq:qtaq3}
    \begin{bmatrix}
    \mathcal{Q}_f^T & & \\ & \mathcal{Q}_u^T & \\ & & \mathcal{Q}_{\lambda}^T \end{bmatrix}
    \begin{bmatrix}
    2 \beta \mathcal{M} & 0 & -\mathcal{M}\\
    0 & \mathcal{M} & \mathcal{K}(\mu)^T\\
    -\mathcal{M} & \mathcal{K}(\mu) & 0
    \end{bmatrix}
    \begin{bmatrix}
    \mathcal{Q}_f & & \\ & \mathcal{Q}_u & \\ & & \mathcal{Q}_{\lambda}\end{bmatrix}=
    \begin{bmatrix}
    \mathcal{Q}_f^T (2 \beta \mathcal{M}) \mathcal{Q}_f &
    0 &
    - \mathcal{Q}_f^T \mathcal{M} \mathcal{Q}_{\lambda}\\
    0 &
    \mathcal{Q}_u^T \mathcal{M} \mathcal{Q}_u &
    \mathcal{Q}_u^T \mathcal{K}(\mu)^T \mathcal{Q}_{\lambda}\\
    -\mathcal{Q}_{\lambda}^T \mathcal{M} \mathcal{Q}_f &
    \mathcal{Q}_{\lambda}^T \mathcal{K}(\mu) \mathcal{Q}_u &
    0
    \end{bmatrix}
    \end{align}}
 \begin{align}\label{eq:qtaq2}\begin{bmatrix}
    \mathcal{Q}_{\bar{x}}^T & \\ & & \mathcal{Q}_{\lambda}^T \end{bmatrix}
    \begin{bmatrix}
     \mathcal{A} & \mathcal{B}(\mu)^T\\
     \mathcal{B}(\mu) &  0
    \end{bmatrix}
    \begin{bmatrix}
    \mathcal{Q}_{\bar{x}} & \\ & & \mathcal{Q}_{\lambda} \end{bmatrix}=
    \begin{bmatrix}
    \mathcal{Q}_{\bar{x}}^T 
    \mathcal{A}
    \mathcal{Q}_{\bar{x}} &
    \mathcal{Q}_{\bar{x}}^T
    \mathcal{B}(\mu)^T 
    \mathcal{Q}_{\lambda}\\
    \mathcal{Q}_{\lambda}^T
    \mathcal{B}(\mu)
    \mathcal{Q}_{\bar{x}} &
    0
    \end{bmatrix}\end{align}
In the online step, simulations are carried out and the matrix $\mathcal{Q}$ projects the full problem onto the reduced space using one of the two methods used to specify the reduced system, Galerkin or Petrov-Galerkin projection  (\ref{eq:galsys})-(\ref{eq:petgalsys}).

\section{Stabilization}\label{sec:stabilization}

The discussion above describes the construction of reduced spaces $X_N=F_N \times U_N$ and $Q_N$, which (so far) give rise to matrices  $\mathcal{Q}_u = [\mathbf{u}(\mu^1),..., \mathbf{u}(\mu^N)]$, and $\mathcal{Q}_f $, $\mathcal{Q}_{\bar{x}}$ and $\mathcal{Q}_{\lambda}$ defined analogously. A question remains about inf-sup stability of the saddle point systems (\ref{eq:qtaq2}) and (\ref{eq:qtaq3}), i.e., whether \begin{align}\label{eq:infsuprewritten}
    \text{inf}_{q \in Q_N} ~ \text{sup}_{\bar{w} \in X_N} ~ \frac{\mathrm{B}(\bar{w}, q; \mu)}{||\bar{w}||_X ||q||_Q} = \min_{\mathbf{q}} \max_{\bar{\mathbf{w}}} \frac{\langle \mathcal{B}(\mu) \bar{\mathbf{w}}, \mathbf{q} \rangle}{||\bar{\mathbf{w}}||_X ||\mathbf{q}||_Q} \geq \beta_0
\end{align}
Note that $X_N =$span$\{\bar{x}(\mu^{(k)}) \}$ where, for homogeneous Dirichlet boundary conditions, $\bar{\mathbf{x}}(\mu^{(k)})= \footnotesize{\begin{bmatrix}
    \mathbf{f}_{\mu^{(k)}} \\ \mathbf{u}_{\mu^{(k)}}
\end{bmatrix}}$ satisfies the weak constraint (\ref{eq:discconstraintstokes}) with $\mathbf{d}_{\mu}=0$. Consequently, the numerator in (\ref{eq:infsuprewritten}) is $0$ for all $\bar{w} \in X_N$ and the reduced problem is not inf-sup stable. We now discuss two techniques designed to address this by enriching $X_N$ and/or $Q_N$.

\subsection{Stabilization by Supremizer}\label{sec:stabsup}

Let $\mathcal{Y}=[\bar{\mathbf{x}}(\mu^{(1)}),...,\bar{\mathbf{x}}(\mu^{(N)})]$ and $\mathcal{L}=[\boldsymbol{\lambda}(\mu^{(1)}),...,\boldsymbol{\lambda}(\mu^{(N)})]$,  where $\{\bar{\mathbf{x}}(\mu^{(j)})\}$ and $\{\boldsymbol{\lambda}(\mu^{(j)})\}$ are the basis vectors in $V_N=X_N \times Q_N$ constructed by Algorithm \ref{alg:greedy}. Following \cite{ballarin2015supremizer}, we describe two ways to construct supremizers to enrich the space determined by $\mathcal{Y}$. The first  of these,  called an ``exact  supremizer" in \cite{ballarin2015supremizer}, produces what is needed but not in  a practical way.  Let $\mu$ be a parameter arising in an online simulation,  and let $\mathbf{r}_j$ be the solution of $\mathcal{A} \mathbf{r}_j= \mathcal{B}(\mu)^T \boldsymbol{\lambda}_j$, $j=1,...,N$. 
(It can be shown that 
$\mathbf{r}_j= \arg \max_{\bar{\mathbf{w}}} \frac{\langle \mathcal{B}(\mu) \bar{\mathbf{w}}, \boldsymbol{\lambda}_j \rangle}{||\bar{\mathbf{w}}||_X}$, whence the name \textit{supremizer}.)
Let $\mathcal{R}=[\mathbf{r}_1,...,\mathbf{r}_N]$, and let the enriched reduced state space be defined as the span of $[\mathcal{Y}, \mathcal{R}]$. 
The operators satisfy $\mathcal{A}\mathcal{R}=\mathcal{B}(\mu)^T\mathcal{L}$, and any member of the enriched space has the form $\bar{\mathbf{w}}=\mathcal{Y} \xi + \mathcal{R} \omega$. 
Therefore, on the enriched reduced spaces, we seek a lower bound on
\begin{align*}
    \min_{\mathbf{q}\in \text{range} (\mathcal{L})} \max_{\bar{\mathbf{w}} \in \text{range} [\mathcal{Y}, \mathcal{R}]} \frac{\langle \mathcal{B}(\mu)\bar{\mathbf{w}}, \mathbf{q} \rangle}{||\bar{\mathbf{w}}||_X ||\mathbf{q}||_Q} = \min_{\alpha} \max_{\xi, \omega} \frac{\langle \mathcal{B}(\mu) [\mathcal{Y} \xi + \mathcal{R} \omega], \mathcal{L} \alpha \rangle}{\langle \mathcal{A} [\mathcal{Y} \xi + \mathcal{R} \omega], [\mathcal{Y} \xi + \mathcal{R} \omega] \rangle^{1/2} \langle \mathcal{L} \alpha, \mathcal{L} \alpha \rangle^{1/2}}
\end{align*}
where $\mathbf{q}=\mathcal{L} \alpha$.  But 
\begin{align*}
    \max_{\xi, \omega} \frac{\langle \mathcal{B}(\mu) [\mathcal{Y} \xi + \mathcal{R} \omega], \mathcal{L} \alpha \rangle}{\langle \mathcal{A} [\mathcal{Y} \xi + \mathcal{R} \omega], [\mathcal{Y} \xi + \mathcal{R} \omega] \rangle^{\frac{1}{2}}}  &\geq \max_{\omega} \frac{\langle \mathcal{B}(\mu) \mathcal{R} \omega, \mathcal{L} \alpha \rangle}{\langle \mathcal{A} \mathcal{R} \omega, \mathcal{R} \omega \rangle^{1/2}} ~~~~~(\xi=0)\\
    &= 
    \max_{\omega} \frac{\langle \omega, \mathcal{R}^T \mathcal{B}(\mu)^T \mathcal{L} \alpha \rangle}{\langle (\mathcal{R}^T \mathcal{A} \mathcal{R})^{1/2} \omega,  (\mathcal{R}^T \mathcal{A} \mathcal{R})^{1/2} \omega \rangle^{1/2}}\\
    &= \max_{\theta}  \frac{\langle \theta, (\mathcal{R}^T \mathcal{A} \mathcal{R})^{-1/2} \mathcal{R}^T \mathcal{B}(\mu)^T \mathcal{L} \alpha \rangle}{||\theta||}~~~~~~~ \text{ where }  ~~ \theta=(\mathcal{R}^T \mathcal{A} \mathcal{R})^{-1/2} \omega\\
    &=||(\mathcal{R}^T \mathcal{A} \mathcal{R})^{-1/2} R^T \mathcal{B}(\mu)^T \mathcal{L} \alpha||
    = \langle (\mathcal{R}^T \mathcal{A} \mathcal{R}) \alpha,   \alpha \rangle^{1/2},
\end{align*}
where the last inequality follows from the fact that $\mathcal{B}(\mu)^T \mathcal{L}= \mathcal{A} \mathcal{R}$.
Thus,  
\begin{align*}
    \max_{\bar{\mathbf{w}} = \mathcal{Y} \xi + \mathcal{R} \omega} \frac{\langle  \mathcal{B}(\mu) \bar{\mathbf{w}}, \mathbf{q} \rangle}{||\bar{\mathbf{w}}||_X ||\mathbf{q}||_Q }& \geq \frac{\langle (\mathcal{R}^T \mathcal{A} \mathcal{R}) \alpha,   \alpha \rangle^{1/2}}{\langle \mathcal{K} \mathcal{L} \alpha, \mathcal{L} \alpha \rangle^{1/2}}= \Big[ \frac{\langle (\mathcal{R}^T \mathcal{A} \mathcal{R}) \alpha,   \alpha \rangle^{1/2}}{\langle \mathcal{A}^{-1} \mathcal{B}(\mu)^T \mathcal{L} \alpha, \mathcal{B}(\mu)^T \mathcal{L} \alpha \rangle^{1/2}} \Big]  \Big[ \frac{\langle \mathcal{A}^{-1} \mathcal{B}(\mu)^T \mathcal{L} \alpha, \mathcal{B}(\mu)^T \mathcal{L} \alpha \rangle^{1/2}}{\langle \mathcal{K} \mathcal{L} \alpha, \mathcal{L} \alpha \rangle^{1/2}}\Big].
\end{align*}
The first of these factors is 
$$\frac{\langle (\mathcal{R}^T \mathcal{A} \mathcal{R}) \alpha,   \alpha \rangle^{1/2}}{\langle \mathcal{A}^{-1} \mathcal{B}(\mu)^T \mathcal{L} \alpha, \mathcal{B}(\mu)^T \mathcal{L} \alpha \rangle^{1/2}}= \frac{\langle (\mathcal{R}^T \mathcal{A} \mathcal{R}) \alpha,   \alpha \rangle^{1/2}}{\langle (\mathcal{R}^T \mathcal{A} \mathcal{R}) \alpha,   \alpha \rangle^{1/2}}=1.$$
From inequality (\ref{eq:schurversion}), the second factor is bounded below by $\beta_0$. Thus, this version of the supremizer produces a div-stable reduced basis.

For reduced-basis methods to be practical, the reduced basis should be constructed in the offline step.  The method just described does not meet this requirement, as the supremizers depend on the parameter $\mu$ used in the online simulation and a new enriched reduced basis must be constructed for each new parameter. A practical variant constructs a set of supremizers $\{ \mathbf{r}_j \}$ that satisfy $\mathcal{A} \mathbf{r}_j = \mathcal{B}(\mu^{(j)})^T \boldsymbol{\lambda}_j$,  where $\{\mu^{(j)}\}$ is the set of parameters chosen during the search, Algorithm \ref{alg:greedy}. This computation can be done in the offline step. The resulting quantities satisfy
\begin{align*}
    \mathcal{A} \mathcal{R} = [\mathcal{B}(\mu^{(1)})^T \boldsymbol{\lambda}_1,...,\mathcal{B}(\mu^{(N)})^T \boldsymbol{\lambda}_N],
\end{align*}
but there is no longer a relation of the form $\mathcal{A} \mathcal{R}=B^T \mathcal{L}$. Consequently, the argument above establishing inf-sup stability of the reduced spaces is not applicable. This approach is described in \cite{ballarin2015supremizer} as constructing \textit{approximate supremizers},  and it is the method we explore in experiments.

Let $\mathcal{Q}_{sup}$ denote the block-diagonal reduced basis matrix with augmentation by supremizer. $\mathcal{Q}_{sup}$  has the form in the right side of (\ref{eq:tworedmat}).
As the basis is constructed, at each step that a snapshot is added, $\mathcal{Q}_{\lambda}$ is updated with snapshot $\boldsymbol{\lambda}_{\mu}$ and $\mathcal{Q}_{\bar{x}}$ is updated with snapshot $\bar{\textbf{x}}_{\mu}$ and supremizer $\textbf{r}_{\mu} = \mathcal{A}^{-1} \mathcal{B}(\mu)^T \boldsymbol{\lambda}_{\mu}$. 
Thus, the matrix $\mathcal{Q}_{\lambda}$ from the naive RB spaces is left unaugmented and $\mathcal{Q}_{\bar{x}}$ is augmented so that its range is span $\{ \bar{\mathbf{x}}_{\mu^1}, \mathbf{r}_{\mu^1},..., \bar{\mathbf{x}}_{\mu^N}, \mathbf{r}_{\mu^N} \}$. For computations, we use a Gram-Schmidt process so that both $\mathcal{Q}_{\bar{x}}$ and $\mathcal{Q}_{\lambda}$ are forced to have orthonormal columns. The matrix $\mathcal{Q}_{\bar{x}}$ has twice as many columns and rows as $\mathcal{Q}_{\lambda}$, and the entire matrix $\mathcal{Q}_{sup}$ has $3N$ columns, where $N$ is the number of snapshots used in Algorithm \ref{alg:greedy}.

Stabilizing in this manner has been considered for both the parameterized optimal control problem as well as for stabilizing PDEs themselves. In \cite{chenquarrozza}, stabilization by supremizer was employed for solving the Stokes control problem. In \cite{elmanforstall_nav}, this method of stabilization is used for stabilizing the  Navier-Stokes equations.

\subsection{Stabilization by Aggregation}\label{sec:stabagg}

Another method of enriching the spaces to ensure inf-sup stability of the RB approximation spaces is to augment by aggregation \cite{nrmq}. Because equivalence of $U$ and $Q$ ensured inf-sup stability at the continuous and discrete levels, both $U_N$ and $Q_N$ are enriched so that they are equivalent and are each updated with both $\textbf{u}(\mu)$ and $\boldsymbol{\lambda}(\mu))$ at each step the RB spaces are updated. 
That is,  the updated spaces and matrices determined by Algorithm \ref{alg:greedy} are given by the following rules:
\begin{itemize}
    \item Let $Z_N=$ span $\{ u_h(\mu^{(n)}), \lambda_h(\mu^{(n)}), n=1,...,N\}$
    \item state: $U_N=Z_N$ such that range$(\mathcal{Q}_{u})=$ range$([\mathbf{u}_{\mu^{(1)}}, \boldsymbol{\lambda}_{\mu^{(1)}},..., \mathbf{u}_{\mu^{(N)}}, \boldsymbol{\lambda}_{\mu^{(N)}}])$ 
    \item control: $F_N=$ span $\{ f(\mu^{(n)}), n=1,...,N\}$ such that range$(\mathcal{Q}_{f})=$ range$([\mathbf{f}_{\mu^{(1)}},..., \mathbf{f}_{\mu^{(N)}}])$
    \item adjoint/Lagrange multiplier: $Q_N=Z_N$ such that $\mathcal{Q}_{\lambda}=\mathcal{Q}_{u}$
\end{itemize}
Inf-sup stability is established as follows  \cite{nrmq};
\begin{align*}
    \min_{\mathbf{q} \in Q_N} ~ \max_{(\mathbf{v}, \mathbf{z}) \in F_N \times U_N} &~ \frac{\langle \mathcal{K}(\mu) \mathbf{z}, \mathbf{q} \rangle - \langle \mathcal{M} \mathbf{v}, \mathbf{q} \rangle}{( \langle 2 \beta \mathcal{M} \mathbf{v}, \mathbf{v} \rangle +  \langle \mathcal{M} \mathbf{z}, \mathbf{z} \rangle)^{1/2} ||\mathbf{q}||_Q}
    \geq \min_{\mathbf{q} \in Q_N} ~\max_{(\mathbf{0}, \mathbf{z}) \in F_N \times U_N} ~ \frac{\langle \mathcal{K}(\mu) \mathbf{z}, \mathbf{q} \rangle }{\langle \mathcal{M} \mathbf{z}, \mathbf{z} \rangle^{1/2} ||\mathbf{q}||_Q}\\
    & \geq \min_{\mathbf{q} \in Q_N} ~ \frac{\langle \mathcal{K}(\mu) \mathbf{q}, \mathbf{q} \rangle }{\langle \mathcal{M} \mathbf{q}, \mathbf{q} \rangle^{1/2} ||\mathbf{q}||_Q}
      \geq  \min_{\mathbf{q}} ~ c_{\Omega} ~\frac{\langle \mathcal{K}(\mu) \mathbf{q}, \mathbf{q} \rangle }{||\mathbf{q}||_Q^2} \geq c_{\Omega} ~\tilde{\alpha}_0 ,
    \end{align*}
    where $c_{\Omega}$, $\tilde{\alpha}_0$ come from the Poincare inequality and coercivity of $a(\cdot, \cdot; \mu)$, respectively. Note that this argument depends on the assumption that $Q_N=U_N$.  The resulting reduced basis is fully constructed in the offline step.

Let $\mathcal{Q}_{agg}$ denote the block-diagonal reduced basis matrix with augmentation by aggregation, which has the structure shown in the left side of (\ref{eq:tworedmat}).   
The RB matrix $\mathcal{Q}_{f}$ from the definition of naive RB spaces is left unaugmented.
The RB matrices $\mathcal{Q}_{u}$ and $\mathcal{Q}_{\lambda}$ are both augmented to ensure  $\mathcal{Q}_{u}=\mathcal{Q}_{\lambda}$. 
It follows that $\mathcal{Q}_u$ and $\mathcal{Q}_{\lambda}$ each have twice as many columns as $\mathcal{Q}_f$ and the total number of columns in $\mathcal{Q}_{agg}$ is $5N$. For computations, as for the supremizer, we use a Gram-Schmidt process so that $\mathcal{Q}_{u}$, $\mathcal{Q}_{f}$ and $\mathcal{Q}_{\lambda}$ have orthonormal columns.

Stabilization by aggregation, also referred to as integration, has been considered strictly for using reduced order modeling to solve parameterized optimal control problems. In \cite{ktgv, nrmq}, this method is used for solving the parameterized elliptic optimal control problems, specifically parameterized diffusion and convection-diffusion control problems.  In \cite{zainib}, aggregation is used along with supremizers for reduced velocity and pressure spaces when solving optimal flow control pipeline problems.

\section{Numerical Results}\label{sec:numres}

\begin{figure}[]
    \centering
    \includegraphics[width=.3\linewidth]{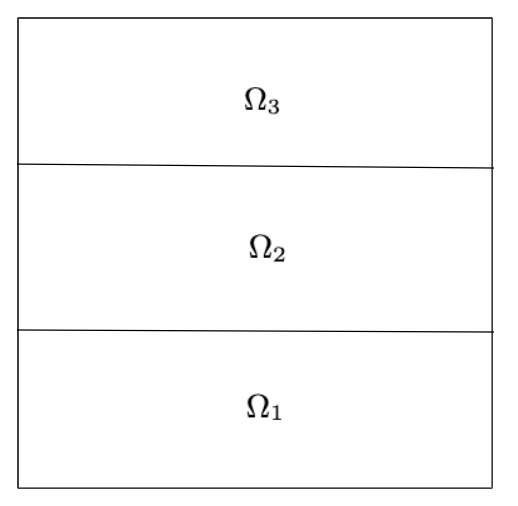}
    \caption{Subdivision of the spatial domain $\Omega$ for the parameterized diffusion control problem  when $N_D=3$. In this case, $\mu \in \mathbb{R}^3$.}
    \label{fig:threestrips}
\end{figure}

In this section, we present experimental results with the two stabilization methods.
We study two benchmark problems, the parametrized diffusion control problem (\ref{eq:pardiffctrl1})-(\ref{eq:pardiffctrl4}) and a parametrized convection-diffusion control problem.
 For the diffusion problem, the spatial domain is $\Omega=[0,1]^2$, with Dirichlet boundary $\Omega_D := [0,1] \times \{ 1 \}$ and Neumann boundary $\Omega_N := \Omega \setminus \Omega_D$. The desired or target state is uniformly equal to one everywhere, $\hat{u}(x, \mu)=1$, and $g(x)=0$. Thus, the target is inconsistent with the Dirichlet boundary condition. The state $u(x, \mu)$  cannot match the target on the Dirichlet part of the boundary and the control $f(x, \mu)$ will require more energy to produce a good approximation of the target state. The domain $D$ is divided into $N_D$ equal-sized horizontal subdomains, where the stratified domain consists of  $N_D$ horizontal strips, as shown in Figure \ref{fig:threestrips} for $N_D=3$.
The diffusion coefficient is piecewise constant on each subdomain
$\sigma(x,\mu) |_{\Omega_k} = \mu_k, k=1:N_D,$
where the parameters $\mu = [\mu_1,..., \mu_{N_D}]^T$ are taken to be independently and uniformly distributed random variables  in $\Gamma := [0.01,1]^{N_D}$. 

The training set (of size $N_{max}$) is used in Algorithm \ref{alg:greedy} to build snapshots that form the basis. The greedy search is implemented with error indicator 
$$ \eta(\mathbf{v}_r)=\frac{|| \mathcal{G}(\mu) \mathbf{v}_r  - \mathbf{b} ||_2}{|| \mathbf{b}||_2}$$
(where $\mathbf{v}_r$ refers to any element in $V_h$) and is carried out until $\eta(\mathbf{v}_r(\mu))$ is less than some desired tolerance for all $\mu \in T$. The basis is tested in the online stage using a set of parameters not in the training set $T$. In all experiments, we used $N_{max}=2000$ randomly chosen training points, regularization parameter $\beta=10^{-2}$ in (\ref{eq:pardiffctrl1}), and spatial discretization consisting of piecewise bilinear finite elements on a uniform grid with $(2^{nc} +1) \times (2^{nc} + 1)$ elements. Computations were done on a Macbook Pro with an Intel 2.2 GHz i7 processor and 16 GB RAM, using MATLAB R2022a, or on a Dell Precision 7820 Tower with an Intel Xeon Silver 1102.1 GHz processor and 64 GB RAM, using MATLAB R2018b. Finite element matrices for the full models were generated using IFISS 3.6 \cite{ifiss}.

\begin{figure}[]
    \centering
    \begin{minipage}{0.48\textwidth}
        \centering
        \includegraphics[width=\linewidth]{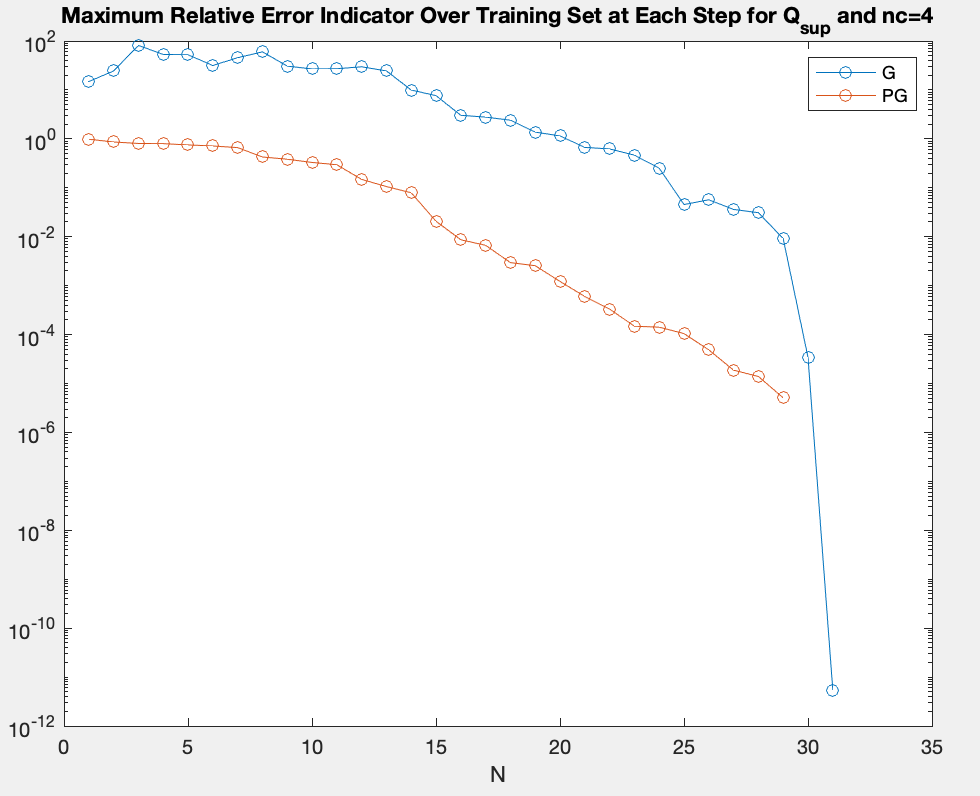}
        \caption{Maximum relative error indicator over the training set as $\mathcal{Q}_{sup}$ is being built with both a Galerkin (G) and Petrov-Galerkin (PG) solve for the parameterized diffusion control problem. The number of basis vectors is $N$. Here, the stopping tolerance for the greedy algorithm is $10^{-7}$, the spatial discretization has $(2^{nc}+1)\times (2^{nc} +1)$ elements where $nc=4$, and the number of subdomains $\Omega_k \subset \Omega$ is $N_D=3$. }
        \label{fig:mr_sup_nc4}
    \end{minipage}%
    \hspace{0.5cm}
     \begin{minipage}{0.48\textwidth}
        \centering
        \includegraphics[width=\linewidth]{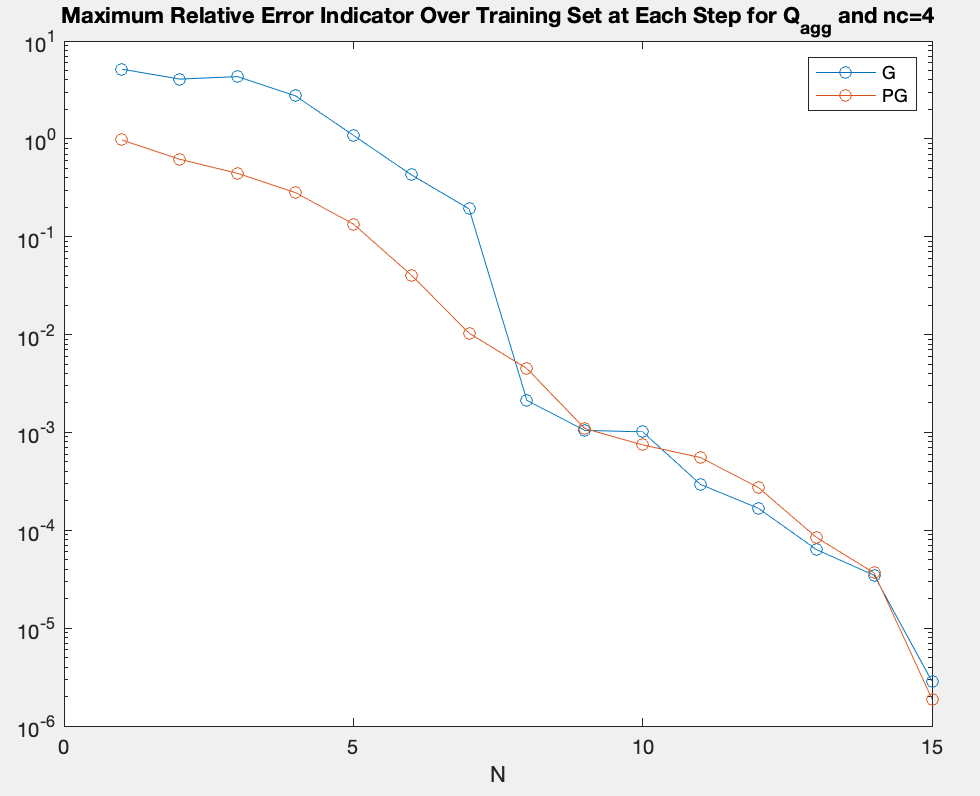}
        \caption{Maximum relative error indicator over the training set as $\mathcal{Q}_{agg}$ is being built with both a Galerkin (G) and Petrov-Galerkin (PG) solve for the parameterized diffusion control problem. The number of basis vectors is $N$. Here, the stopping tolerance for the greedy algorithm is $10^{-7}$, the spatial discretization has $(2^{nc}+1)\times (2^{nc} +1)$ elements where $nc=4$, and the number of subdomains $\Omega_k \subset \Omega$ is $N_D=3$.}
        \label{fig:mr_agg_nc4}
    \end{minipage}
\end{figure}

\begin{figure}[]
    \centering
    \begin{minipage}{0.48\textwidth}
        \centering
        \includegraphics[width=\linewidth]{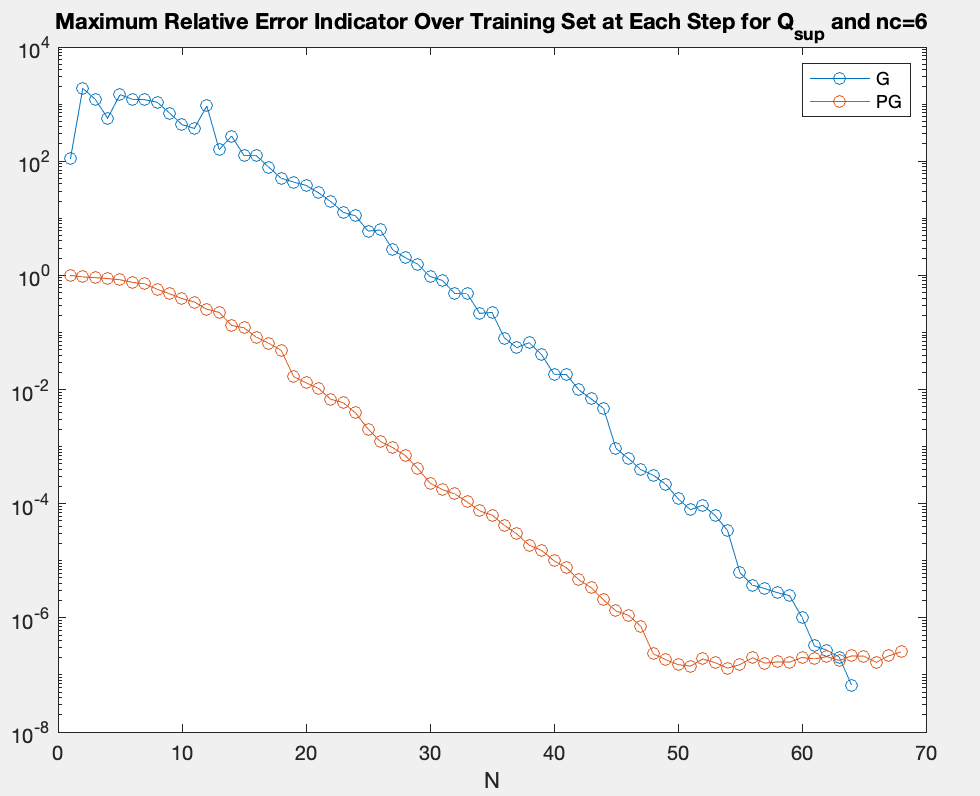}
        \caption{Maximum relative error indicator over the training set as $\mathcal{Q}_{sup}$ is being built with both a Galerkin (G) and Petrov-Galerkin (PG) solve for the parameterized diffusion control problem. The number of basis vectors is $N$. Here, the stopping tolerance for the greedy algorithm is $10^{-7}$, the spatial discretization has $(2^{nc}+1)\times (2^{nc} +1)$ elements where $nc=6$, and the number of subdomains $\Omega_k \subset \Omega$ is $N_D=3$. The PG reduction fails to converge.}
        \label{fig:mr_sup_nc6}
    \end{minipage}%
    \hspace{0.5cm}
     \begin{minipage}{0.48\textwidth}
        \centering
        \includegraphics[width=\linewidth]{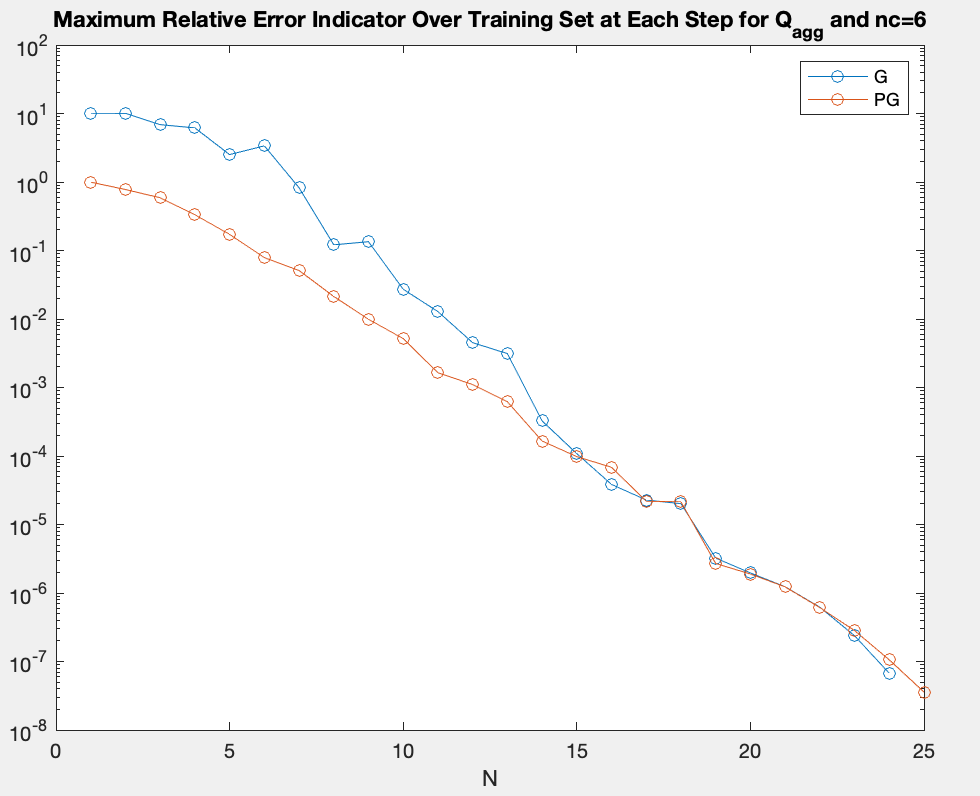}
        \caption{Maximum relative error indicator over the training set as $\mathcal{Q}_{agg}$ is being built with both a Galerkin (G) and Petrov-Galerkin (PG) solve for the parameterized diffusion control problem. The number of basis vectors is $N$. Here, the stopping tolerance for the greedy algorithm is $10^{-7}$, the spatial discretization has $(2^{nc}+1)\times (2^{nc} +1)$ elements where $nc=6$, and the number of subdomains $\Omega_k \subset \Omega$ is $N_D=3$.\vspace{3mm}}
        \label{fig:mr_agg_nc6}
    \end{minipage}
\end{figure}

We begin by examining how training, i.e., Algorithm \ref{alg:greedy}, behaves for the two methods of stabilization, using both Galerkin and Petrov-Galerkin formulations of the reduced problem. Figure \ref{fig:mr_sup_nc4} shows the maximum relative error indicator for the supremizer over the training set $T$ as $\mathcal{Q}_{sup}$ is constructed, for $nc=4$ and $N_D=3$. The tolerance for the greedy algorithm is $10^{-7}$. Figure \ref{fig:mr_agg_nc4} shows the analogous results for aggregation. Note that the Petrov-Galerkin formulation (\ref{eq:petgalsys}) corresponds to the normal equations associated with the (residual) error indicator $\eta(\mu)$, so  that as the reduced basis grows, the error indicator monotonically decreases. This is not true for Galerkin formulation. It can be seen from Figure \ref{fig:mr_sup_nc4} that with the given error indicator, the Petrov-Galerkin method is more effective than the Galerkin method when the supremizer is used for stabilization. In contrast, the behavior of the two reduced models is closer when aggregation is used. However, the behaviors of the two reduction strategies are very similar (and are close to identical for certain values of $N$) as the maximum error indicator values approach the desired tolerance.  Similar results are shown in Figures \ref{fig:mr_sup_nc6} and \ref{fig:mr_agg_nc6} for a finer spatial mesh. For the supremizer, the Petrov-Galerkin formulation produces better maximum relative error indicator values for most values of $N$. However, there is one significant difference between Figures \ref{fig:mr_sup_nc4} and \ref{fig:mr_sup_nc6}. On the finer mesh (Figure \ref{fig:mr_sup_nc6}), the Petrov-Galerkin formulation fails to meet the tolerance (of $10^{-7}$), whereas the Galerkin formulation succeeds. For aggregation, there is better behavior of the Petrov-Galerkin formulation initially and very similar behavior of the two formulations as the tolerance is reached. Both formulations are shown to be convergent in Figure \ref{fig:mr_agg_nc6}, although we note that on a finer mesh ($nc=7$), the Petrov-Galerkin formulation also failed.

\begin{figure}[]
    \centering
    \begin{minipage}{0.48\textwidth}
        \centering
        \includegraphics[width=\linewidth]{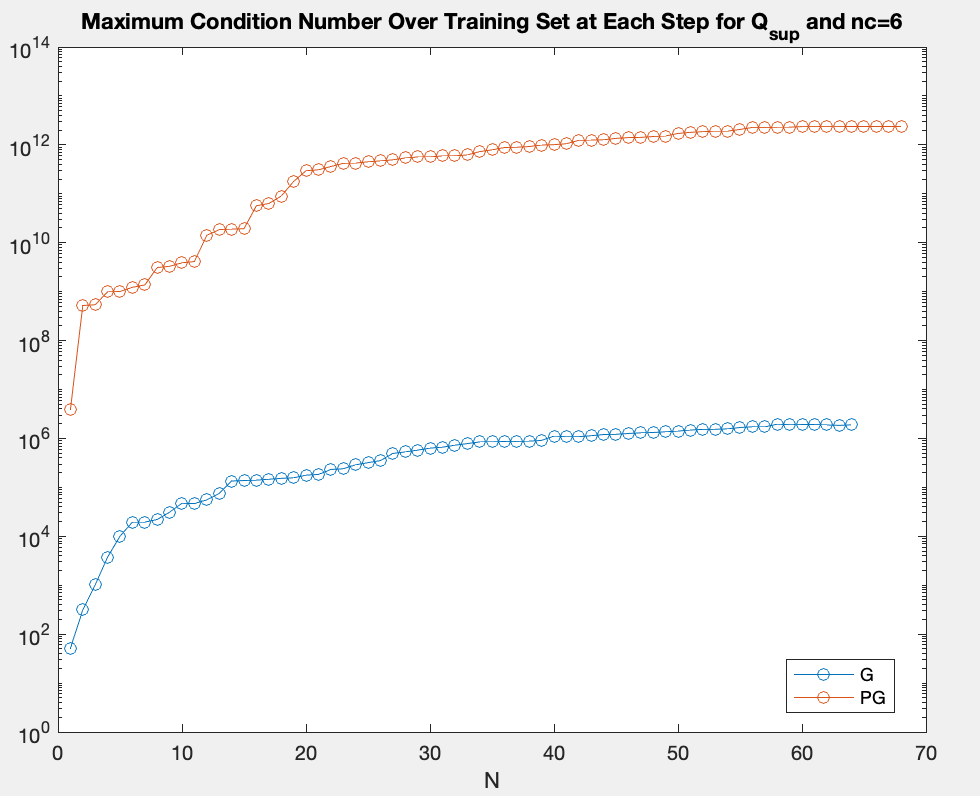}
        \caption{Maximum condition number of the reduced system $\mathcal{Q}^T \mathcal{G}(\mu) \mathcal{Q}$ in the Galerkin (G) case and $(\mathcal{G}(\mu) \mathcal{Q})^{T} (\mathcal{G}(\mu) \mathcal{Q})$ in the Petrov-Galerkin (PG) case over all parameters in the training set as $\mathcal{Q}_{sup}$ is being built for the parameterized diffusion control problem. The number of basis vectors is $N$. Here, the stopping tolerance for the greedy algorithm is $10^{-7}$, the spatial discretization has $(2^{nc}+1)\times (2^{nc} +1)$ elements where $nc=6$, and the number of subdomains $\Omega_k \subset \Omega$ is $N_D=3$.} 
        \label{fig:cond_sup_nc6}
    \end{minipage}%
    \hspace{0.5cm}
     \begin{minipage}{0.48\textwidth}
        \centering
        \includegraphics[width=\linewidth]{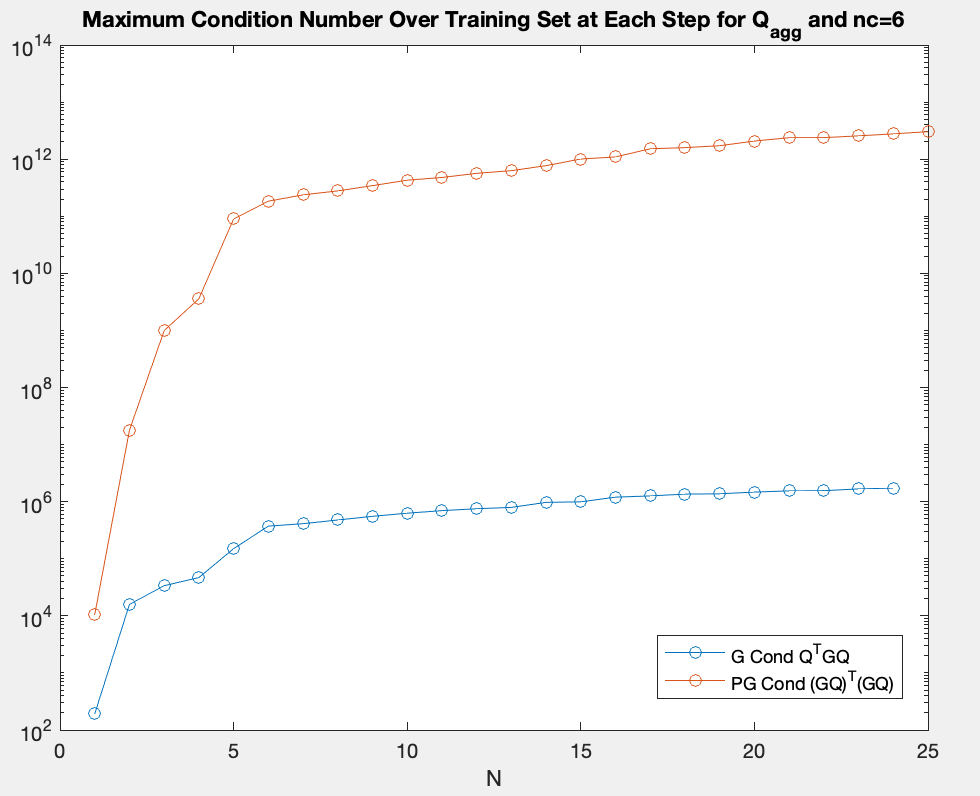}
        \caption{Maximum condition number of the reduced system $\mathcal{Q}^T \mathcal{G}(\mu) \mathcal{Q}$ in the Galerkin (G) case and $(\mathcal{G}(\mu) \mathcal{Q})^{T} (\mathcal{G}(\mu) \mathcal{Q})$ in the Petrov-Galerkin (PG) case over all parameters in the training set as $\mathcal{Q}_{agg}$ is built for the parameterized diffusion control problem. The number of basis vectors is $N$. Here, the stopping tolerance for the greedy algorithm is $10^{-7}$, the spatial discretization has $(2^{nc}+1)\times (2^{nc} +1)$ elements where $nc=6$, and the number of subdomains $\Omega_k \subset \Omega$ is $N_D=3$. \vspace{4mm}}
        \label{fig:cond_agg_nc6}
    \end{minipage}
\end{figure}

This brings about a question of robustness of the Petrov-Galerkin formulation. Figures \ref{fig:cond_sup_nc6} and \ref{fig:cond_agg_nc6} show the maximum condition numbers of the reduced matrices $\mathcal{Q}^T \mathcal{G}(\mu) \mathcal{Q}$ (Galerkin) and $(\mathcal{G}(\mu) \mathcal{Q})^T(\mathcal{G}(\mu) \mathcal{Q})$ (Petrov-Galerkin) over the training set $T$. The latter matrices become severely ill-conditioned as the search proceeds, and we attribute its failure to produce a suitable reduced basis to ill-conditioning. Note that this is a common issue in linear algebra: a matrix $\mathcal{A}^T \mathcal{A}$ associated with the normal equations has condition number equal to the square of the condition number of $\mathcal{A}$. The trends seen in Figures  \ref{fig:cond_sup_nc6} and \ref{fig:cond_agg_nc6}, for the two stabilization methods, are consistent as the basis size varies. In Table \ref{tab:galpetgalcond}, results for the maximum condition number of the reduced systems over the training set for all $N$ are given for both formulations for various spatial meshes. Both the supremizer and aggregation lead to lack of convergence of the greedy algorithm for the Petrov-Galerkin formulation when $nc=7$. The Galerkin formulation never fails to converge.

\begin{table}[h]
\centering
\begin{small}
\begin{tabular}{|c|cc|cc|}
\hline
     & \multicolumn{2}{c|}{Galerkin}              & \multicolumn{2}{c|}{Petrov-Galerkin}       \\ \hline
$nc$ & \multicolumn{1}{c|}{$\mathcal{Q}_{sup}$} & $\mathcal{Q}_{agg}$ & \multicolumn{1}{c|}{$\mathcal{Q}_{sup}$} & $\mathcal{Q}_{agg}$ \\ \hline
3    & \multicolumn{1}{c|}{$4.7e\!\!+\!\!04$} & $4.7e\!\!+\!\!04$ & \multicolumn{1}{c|}{$2.3e\!\!+\!\!09$} & $2.2e\!\!+\!\!09$ \\ \hline
4    & \multicolumn{1}{c|}{$2.0e\!\!+\!\!05$} & $2.0e\!\!+\!\!05$ & \multicolumn{1}{c|}{$2.2e\!\!+\!\!10$} & $3.8e\!\!+\!\!10$ \\ \hline
5    & \multicolumn{1}{c|}{$6.0e\!\!+\!\!05$} & $5.5e\!\!+\!\!05$ & \multicolumn{1}{c|}{$2.5e\!\!+\!\!11$} & $3.0e\!\!+\!\!11$ \\ \hline
6    & \multicolumn{1}{c|}{$2.0e\!\!+\!\!06$} & $1.7e\!\!+\!\!06$ & \multicolumn{1}{c|}{-}         & $3.0e\!\!+\!\!12$ \\ \hline
7    & \multicolumn{1}{c|}{$6.1e\!\!+\!\!06$} & $5.0e\!\!+\!\!06$ & \multicolumn{1}{c|}{-}         & -         \\ \hline
\end{tabular}
\end{small}
    \caption{Maximum condition number of the reduced system $\mathcal{Q}^T \mathcal{G}(\mu) \mathcal{Q}$ in the Galerkin case and $(\mathcal{G}(\mu) \mathcal{Q})^{T} (\mathcal{G}(\mu)\mathcal{Q})$ in the Petrov-Galerkin case over all parameters in the training set and all steps of the greedy algorithm for the parameterized diffusion control problem as $nc$ is refined, where the spatial mesh has $(2^{nc} +1) \times (2^{nc} +1) $ discrete elements, tolerance is $10^{-7}$, $N_{max}=2000$ and $N_D=3$. Empty cells represent lack of convergence to desired tolerance.}
\label{tab:galpetgalcond}
\end{table}

\afterpage{
\begin{table}[h]
\centering
\begin{small}
\begin{tabular}{|c|cccccc|cccccc|}
\hline
     & \multicolumn{6}{c|}{Galerkin Projection}                                                                                                                                            & \multicolumn{6}{c|}{Petrov-Galerkin Projection}                                                                                                                                     \\ \hline
     & \multicolumn{2}{c|}{$N$}                                        & \multicolumn{2}{c|}{Columns}                                    & \multicolumn{2}{c|}{\begin{tabular}[c]{@{}c@{}}Maximum Relative\\ Error Indicator\end{tabular}}        & \multicolumn{2}{c|}{$N$}                                        & \multicolumn{2}{c|}{Columns}                                    & \multicolumn{2}{c|}{\begin{tabular}[c]{@{}c@{}}Maximum Relative\\ Error Indicator\end{tabular}}        \\ \hline
$nc$ & \multicolumn{1}{c|}{$Q_{sup}$} & \multicolumn{1}{c|}{$Q_{agg}$} & \multicolumn{1}{c|}{$Q_{sup}$} & \multicolumn{1}{c|}{$Q_{agg}$} & \multicolumn{1}{c|}{$Q_{sup}$} & $Q_{agg}$ & \multicolumn{1}{c|}{$Q_{sup}$} & \multicolumn{1}{c|}{$Q_{agg}$} & \multicolumn{1}{c|}{$Q_{sup}$} & \multicolumn{1}{c|}{$Q_{agg}$} & \multicolumn{1}{c|}{$Q_{sup}$} & $Q_{agg}$ \\ \hline
3    & \multicolumn{1}{c|}{15}        & \multicolumn{1}{c|}{8}         & \multicolumn{1}{c|}{45}        & \multicolumn{1}{c|}{40}        & \multicolumn{1}{c|}{$4.4e\!\!-\!\!14$} & $4.2e\!\!-\!\!14$ & \multicolumn{1}{c|}{15}        & \multicolumn{1}{c|}{8}         & \multicolumn{1}{c|}{45}        & \multicolumn{1}{c|}{40}        & \multicolumn{1}{c|}{$1.8e\!\!-\!\!10$} & $4.3e\!\!-\!\!12$ \\ \hline
4    & \multicolumn{1}{c|}{31}        & \multicolumn{1}{c|}{16}        & \multicolumn{1}{c|}{93}        & \multicolumn{1}{c|}{80}        & \multicolumn{1}{c|}{$6.6e\!\!-\!\!12$} & $2.0e\!\!-\!\!13$ & \multicolumn{1}{c|}{30}        & \multicolumn{1}{c|}{16}        & \multicolumn{1}{c|}{90}        & \multicolumn{1}{c|}{80}        & \multicolumn{1}{c|}{$1.1e\!\!-\!\!07$} & $4.0e\!\!-\!\!11$ \\ \hline
5    & \multicolumn{1}{c|}{57}        & \multicolumn{1}{c|}{23}        & \multicolumn{1}{c|}{171}       & \multicolumn{1}{c|}{115}       & \multicolumn{1}{c|}{$8.7e\!\!-\!\!08$} & $6.9e\!\!-\!\!08$ & \multicolumn{1}{c|}{48}        & \multicolumn{1}{c|}{23}        & \multicolumn{1}{c|}{144}       & \multicolumn{1}{c|}{115}       & \multicolumn{1}{c|}{$3.1e\!\!-\!\!08$} & $7.0e\!\!-\!\!08$ \\ \hline
6    & \multicolumn{1}{c|}{64}        & \multicolumn{1}{c|}{24}        & \multicolumn{1}{c|}{192}       & \multicolumn{1}{c|}{120}       & \multicolumn{1}{c|}{$3.5e\!\!-\!\!07$} & $4.8e\!\!-\!\!08$ & \multicolumn{1}{c|}{-}         & \multicolumn{1}{c|}{25}        & \multicolumn{1}{c|}{-}         & \multicolumn{1}{c|}{125}       & \multicolumn{1}{c|}{-}         & $4.2e\!\!-\!\!08$ \\ \hline
7    & \multicolumn{1}{c|}{67}        & \multicolumn{1}{c|}{25}        & \multicolumn{1}{c|}{201}       & \multicolumn{1}{c|}{125}       & \multicolumn{1}{c|}{$2.8e\!\!-\!\!07$} & $4.8e\!\!-\!\!08$ & \multicolumn{1}{c|}{-}         & \multicolumn{1}{c|}{-}         & \multicolumn{1}{c|}{-}         & \multicolumn{1}{c|}{-}         & \multicolumn{1}{c|}{-}         & -         \\ \hline
\end{tabular}
\end{small}
\caption{Comparison of number of snapshots $N$, number of columns in the reduced basis, and maximum relative error indicator over the verification set for $\mathcal{Q}_{sup}$ and $\mathcal{Q}_{agg}$ with Galerkin and Petrov-Galerkin formulations for the parameterized diffusion control problem. Here, the verification set has 500 parameters, the domain $\Omega$ has $N_D=3$ subdomains, the spatial mesh has $(2^{nc} +1) \times (2^{nc} +1) $ discrete elements, $N_{max}=2000$ and the stopping tolerance for the greedy algorithm is $10^{-7}$. Empty cells correspond to cases where the greedy search failed to reach this tolerance.}
\label{tab:galcolmaxrelerr}
\end{table}
}

Now we examine the performance of the reduced bases in the online stage. In Table \ref{tab:galcolmaxrelerr}, the accuracy of the reduced solution tested over $500$ parameters not in the training set $T$ is shown for a variety of meshes with fixed $N_D=3$ and both Galerkin and Petrov-Galerkin formulations. 
Similar results are given for $N_D=10$ in Table \ref{tab:nd10results} where only results for the Galerkin formulation are shown.  We summarize the trends observed in these experiments, as follows: 
\begin{itemize} \itemsep -2pt
    \item The number of snapshots, as well as the size of the reduced basis, is smaller for $\mathcal{Q}_{agg}$ than for $\mathcal{Q}_{sup}$.
    \item The Petrov-Galerkin formulation is less robust than Galerkin formulation, in the sense that the greedy search failed to reach the stopping tolerance for finer spatial meshes. 
    \item The sizes of the reduced bases tend toward asymptotic limits as the spatial mesh is refined. For example, this limit is approximately $N=25$ snapshots (reduced matrix of order 125) for $\mathcal{Q}_{agg}$ with Galerkin search.
    \item The sizes of the reduced bases are larger for the larger number of parameters $N_{D}=10$.  (Compare Tables \ref{tab:galcolmaxrelerr} and \ref{tab:nd10results}.\footnote{For certain tests in the online stage, shown in these tables, the relative error indicator values produced over the testing set of parameters were slightly greater than the prescribed tolerance. This can be explained by the size of the training set chosen. Increasing the number of training parameters to 3000 resulted in maximum relative error indicator values below the tolerance in all cases.})
\end{itemize}

\begin{figure}[H]
    \centering
    \includegraphics[width=.3\linewidth]{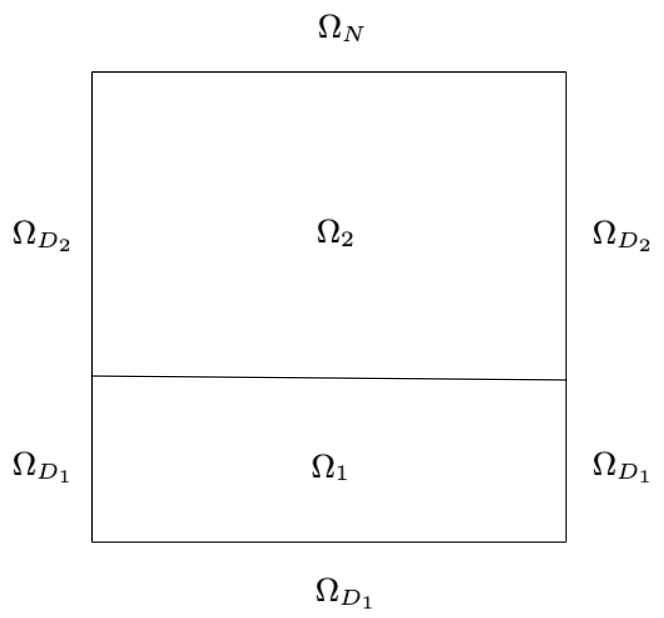}
    \caption{Picture of the spatial domain $\Omega$ and subdomains $\Omega_1$ and $\Omega_2$ for the convection-diffusion control problem.}
    \label{fig:cdspatial}
\end{figure}

For the second benchmark problem,   we consider a variant of the Graetz convection-diffusion problem presented in \cite{nrmq}: find state $u$ and control $f$ such that
\begin{align*}
\min_{u,f} ~ \frac{1}{2} \norm{u(x, \mu) - \hat{u}(x, \mu)}_{L_{2}(\Omega)}^{2} + \frac{\beta}{2} \norm{f(x,\mu)}_{L_{2}(\Omega)}^{2}\\
\text{subject to} \quad -\mu_1 \bigtriangleup  u(x, \mu)  +  \textbf{w} \cdot \bigtriangledown u(x, \mu)~ = ~ f(x, \mu) ~ \text{ in }  ~ \Omega \times \Gamma, \\
\quad \quad \text{such that} \quad \quad \quad  u(x, \mu)  = 1 \quad  ~ \text{on} ~ \partial \Omega_{D_1} \times \Gamma. \quad \quad    \quad \quad \quad \quad \\
\quad \quad \quad \quad  u(x, \mu)  = 2 \quad  ~ \text{on} ~ \partial \Omega_{D_2} \times \Gamma. \quad \quad    \quad \quad \quad \quad \\
\quad \quad \quad \quad \quad  \mu_1 \frac{\partial  u(x, \mu)}{\partial n}  = 0 \quad  ~ \text{on} ~ \partial \Omega_{N} \times \Gamma.    \quad \quad \quad \quad \quad \quad
\end{align*}
Here, $\Omega=[0,1]^2 \subset \R^2$ is the spatial domain (shown in Figure \ref{fig:cdspatial}) subdivided into $\Omega_1=[0,1] \times [0,0.3]$, $\Omega_2=[0,1] \times (0.3,1]$. 
The parameter vector $\mu \in  \Gamma:=[\frac{1}{20}, \frac{1}{3}] \times [0.5,1.5] \times [1.5, 2.5] \in \R^3$ is associated with the diffusion coefficient and the desired state $\hat{u}$ such that $\mu_1$ is the diffusion coefficient, $\hat{u}= \mu_2$ in $\Omega_1=[0,1] \times [0,0.3]$ and $\hat{u}=\mu_3$ in $\Omega_2=[0,1] \times [0.3,1]$. Here, $\textbf{w} = [x_2(1-x_2),0]^T$.

\begin{table}[]
\centering
\begin{small}
\begin{tabular}{|c|cc|cc|cc|}
\hline
     & \multicolumn{2}{c|}{N}                                                        & \multicolumn{2}{c|}{Columns}                                                  & \multicolumn{2}{c|}{\begin{tabular}[c]{@{}c@{}}Maximum Relative\\ Error Indicator\end{tabular}}                       \\ \hline
$nc$ & \multicolumn{1}{c|}{$\mathcal{Q}_{sup}$}  & $\mathcal{Q}_{agg}$ & \multicolumn{1}{c|}{$\mathcal{Q}_{sup}$} &  $\mathcal{Q}_{agg}$ & \multicolumn{1}{c|}{$\mathcal{Q}_{sup}$} & $\mathcal{Q}_{agg}$     \\ \hline
3    & \multicolumn{1}{c|}{15}            & 8      & \multicolumn{1}{c|}{45}                   & 40      & \multicolumn{1}{c|}{$3.9e\!\!-\!\!13$ } & $3.9e\!\!-\!\!14$ \\ \hline
4    & \multicolumn{1}{c|}{31}                    & 16      & \multicolumn{1}{c|}{93}                    & 80      & \multicolumn{1}{c|}{$5.4e\!\!-\!\!12$ } &  $2.1e\!\!-\!\!13$ \\ \hline
5    & \multicolumn{1}{c|}{63}                    & 32      & \multicolumn{1}{c|}{189}                 & 160      & \multicolumn{1}{c|}{$1.3e\!\!-\!\!10$}  & $1.0e\!\!-\!\!12$ \\ \hline
6    & \multicolumn{1}{c|}{126}                   & 52     & \multicolumn{1}{c|}{372}                  & 260     & \multicolumn{1}{c|}{$2.2e\!\!-\!\!08$} & $4.6e\!\!-\!\!08$ \\ \hline
7    & \multicolumn{1}{c|}{173}                  &   $56$      & \multicolumn{1}{c|}{519}                  &   280      & \multicolumn{1}{c|}{$4.6e\!\!-\!\!07$}  &    $7.9e\!\!-\!\!08$  \\ \hline
\end{tabular}
\end{small}
\caption{Comparison of columns and maximum relative error indicator over the verification set for $\mathcal{Q}_{sup}$ and $\mathcal{Q}_{agg}$ with Galerkin projection for the parameterized diffusion control problem. Here, the verification set had 500 parameters, the domain $\Omega$ has $N_D=10$ subdomains, the spatial mesh has $(2^{nc} +1) \times (2^{nc} +1) $ discrete elements, $N_{max}=2000$ and the stopping tolerance for the greedy algorithm is $10^{-7}$.}
\label{tab:nd10results}
\end{table}

The performance of the reduced basis in the online stage is reported for the parameterized convection-diffusion benchmark problem. Table \ref{tab:convdiffresults} is the analog of Table \ref{tab:galcolmaxrelerr} for the parameterized diffusion control problem. As before, stabilizing with aggregation produces a smaller reduced basis, making it more efficient than stabilizing with the supremizer in all cases. In all cases, the basis sizes tend to asymptotic limits but the overall size of the bases are smaller than seen for diffusion control. For this example, the Petrov-Galerkin method was somewhat more efficient than the Galerkin method, although both choices produced solutions with reduced bases of small size. 

\begin{table}[H]
\centering
\begin{small}
\begin{tabular}{|c|cccccc|cccccc|}
\hline
     & \multicolumn{6}{c|}{Galerkin Solve}                                                                                                                                            & \multicolumn{6}{c|}{Petrov-Galerkin Solve}                                                                                                                                     \\ \hline
     & \multicolumn{2}{c|}{$N$}                                        & \multicolumn{2}{c|}{Columns}                                    & \multicolumn{2}{c|}{\begin{tabular}[c]{@{}c@{}}Maximum Relative\\ Error Indicator\end{tabular}}        & \multicolumn{2}{c|}{$N$}                                        & \multicolumn{2}{c|}{Columns}                                    & \multicolumn{2}{c|}{\begin{tabular}[c]{@{}c@{}}Maximum Relative\\ Error Indicator\end{tabular}}        \\ \hline
$nc$ & \multicolumn{1}{c|}{$Q_{sup}$} & \multicolumn{1}{c|}{$Q_{agg}$} & \multicolumn{1}{c|}{$Q_{sup}$} & \multicolumn{1}{c|}{$Q_{agg}$} & \multicolumn{1}{c|}{$Q_{sup}$} & $Q_{agg}$ & \multicolumn{1}{c|}{$Q_{sup}$} & \multicolumn{1}{c|}{$Q_{agg}$} & \multicolumn{1}{c|}{$Q_{sup}$} & \multicolumn{1}{c|}{$Q_{agg}$} & \multicolumn{1}{c|}{$Q_{sup}$} & $Q_{agg}$ \\ \hline
3    & \multicolumn{1}{c|}{17}        & \multicolumn{1}{c|}{10}        & \multicolumn{1}{c|}{51}        & \multicolumn{1}{c|}{50}        & \multicolumn{1}{c|}{$2.9e\!\!-\!\!05$} & $3.0e\!\!-\!\!05$ & \multicolumn{1}{c|}{13}        & \multicolumn{1}{c|}{8}         & \multicolumn{1}{c|}{39}        & \multicolumn{1}{c|}{40}        & \multicolumn{1}{c|}{$4.6e\!\!-\!\!05$} & $6.9e\!\!-\!\!05$ \\ \hline
4    & \multicolumn{1}{c|}{16}        & \multicolumn{1}{c|}{9}         & \multicolumn{1}{c|}{48}        & \multicolumn{1}{c|}{45}        & \multicolumn{1}{c|}{$8.4e\!\!-\!\!05$} & $7.4e\!\!-\!\!05$ & \multicolumn{1}{c|}{11}         & \multicolumn{1}{c|}{7}         & \multicolumn{1}{c|}{33}        & \multicolumn{1}{c|}{35}        & \multicolumn{1}{c|}{$8.7e\!\!-\!\!05$} & $7.9e\!\!-\!\!05$ \\ \hline
5    & \multicolumn{1}{c|}{19}        & \multicolumn{1}{c|}{8}         & \multicolumn{1}{c|}{57}        & \multicolumn{1}{c|}{40}        & \multicolumn{1}{c|}{$4.8e\!\!-\!\!05$} & $6.3e\!\!-\!\!05$ & \multicolumn{1}{c|}{9}         & \multicolumn{1}{c|}{5}         & \multicolumn{1}{c|}{27}        & \multicolumn{1}{c|}{25}        & \multicolumn{1}{c|}{$8.8e\!\!-\!\!05$} & $9.1e\!\!-\!\!05$ \\ \hline
6    & \multicolumn{1}{c|}{19}        & \multicolumn{1}{c|}{5}         & \multicolumn{1}{c|}{57}        & \multicolumn{1}{c|}{25}        & \multicolumn{1}{c|}{$5.2e\!\!-\!\!05$} & $9.9e\!\!-\!\!05$ & \multicolumn{1}{c|}{6}         & \multicolumn{1}{c|}{3}         & \multicolumn{1}{c|}{18}        & \multicolumn{1}{c|}{15}        & \multicolumn{1}{c|}{$9.0e\!\!-\!\!05$} & $9.9e\!\!-\!\!05$ \\ \hline
7    & \multicolumn{1}{c|}{20}        & \multicolumn{1}{c|}{5}         & \multicolumn{1}{c|}{60}        & \multicolumn{1}{c|}{25}        & \multicolumn{1}{c|}{$8.1e\!\!-\!\!05$} & $4.1e\!\!-\!\!05$ & \multicolumn{1}{c|}{4}         & \multicolumn{1}{c|}{2}         & \multicolumn{1}{c|}{12}        & \multicolumn{1}{c|}{10}        & \multicolumn{1}{c|}{$8.2e\!\!-\!\!05$} & $8.5\!\!-\!\!05$ \\ \hline
\end{tabular}
\end{small}
\caption{Comparison of columns, number of basis vectors $N$ and maximum relative error indicator over the verification set for $\mathcal{Q}_{sup}$ and $\mathcal{Q}_{agg}$ with Galerkin and Petrov-Galerkin projection for the parameterized convection-diffusion control problem. Here, the verification set had 500 parameters, the spatial mesh has $(2^{nc} +1) \times (2^{nc} +1) $ discrete elements, $N_{max}=2000$ and the stopping tolerance for the greedy algorithm is $10^{-4}$.}
\label{tab:convdiffresults}
\end{table}

\section{Concluding Remarks}\label{sec:conclusion}

Stabilization of reduced basis models to solve optimal control problems can be handled in multiple ways. When these models are implemented with block diagonal reduced basis matrices, enrichment of the reduced basis spaces is required to ensure well-posedness. Two ways of handling this enrichment and stabilizing the reduced bases are stabilization by aggregation and stabilization using the supremizer function. While both are suitable for ensuring stability, augmenting by aggregation is a superior choice in numerous ways. In particular, we showed that for several examples, aggregation leads to smaller reduced bases than the supremizer, and it is also more robust with respect to convergence.

We also note that this study considers these approaches for one class of problems, arising from optimal control with PDE constraints. Reduced basis methods are also useful in other settings, for example, for parametrized versions of models of computational fluid dynamics with an incompressibility constraint, such as the Stokes and Navier-Stokes equations. One drawback to augmentation by aggregation is that it has only been implemented for stabilizing reduced order models for optimal control problems, whereas augmenting using the supremizer function has proven useful for solving PDEs like the Stokes equations as well. We will consider these issues in a follow-on study \cite{davieelman2}.

\bibliography{my_abbrev_bib}

\end{document}